\documentclass[12pt]{article}

\usepackage{amsmath,amssymb}

\begin{document}

\newcommand{\N}{\mathbb{N}}
\newcommand{\R}{\mathbb{R}}
\newcommand{\Q}{\mathbb{Q}}

\newcommand{\uball}{\mathbb{B}}

\newcommand{\dom}{{\rm dom}}
\newcommand{\Dom}{{\rm Dom}}
\newcommand{\rg}{{\rm Rg}}
\newcommand{\Gr}{{\rm Gr}}
\newcommand{\co}{{\rm co}}
\newcommand{\diam}{{\rm diam\,}}
\newcommand{\epi}{{\rm epi}}
\newcommand{\lev}{{\rm Lev}}
\newcommand{\Rk}{{\R^{\bf k}}}
\newcommand{\Lip}{{\rm Lip\,}}

\newcommand{\Rinf}{{\R\cup\{+\infty\}}}
\newcommand{\Int}{{\rm Int}}
\newcommand{\Cl}{{\rm Cl}}
\newcommand{\eps}{{\varepsilon}}
\newcommand{\black}{\rule{6pt}{6pt}\par\medskip\par}

\newcommand{\ds}{\displaystyle}
\newcommand{\la}{\lambda}
\newcommand{\al}{\alpha}
\newcommand{\m}{\mu}
\newcommand{\lan}{\langle}
\newcommand{\ran}{\rangle}
\newcommand{\tos}{\rightrightarrows}

\newtheorem{theo}{Theorem}[section]
\newtheorem{defi}[theo]{Definition}
\newtheorem{rem}[theo]{Remark}
\newtheorem{lemma}[theo]{Lemma}
\newtheorem{prop}[theo]{Proposition}
\newtheorem{ex}[theo]{Example}
\newtheorem{cor}[theo]{Corollary}

\newcommand{\proof}{\noindent{\bf Proof: }}
\newfont{\sss}{msam10}
\newcommand{\eproof}{\hfill {\sss \char4}\par\medskip}

\title{\LARGE\bf A general scheme for separably reducible properties\thanks{Part of this article has been done during visits of the third author in Technion, Haifa,  within the framework of a joint project between the Bulgarian Academy of Sciences and Israeli Academy of Sciences and Humanities}}

\author{\bf M. Fabian, A. Ioffe and J.P. Revalski\thanks{The first author has been supported by grants GA\v CR 20-2230L and RVO 67985840. } }

\date{}

\maketitle

\begin{abstract}
We propose a general scheme for studying separably reducible properties in metric spaces and then apply it to obtain separable determinacy  of Lipschitz property and the separable determinacy of slopes.
\end{abstract}


\medskip


\section{Introduction and  preliminaries}

The aim of this article is to propose a possible general scheme for investigating the separable reducibility of some properties in a metric space $(X,d)$. Let us remind that a property ${\cal P}$ considered at the points  of the metric space $X$ is said to be {\it separably reducible} (or, equivalently, {\it separably determined}) if there is a family ${\cal F}$ of closed separable subsets of $X$ such that $\cup\{Y: Y\in {\cal F}\}=X$ and  for every $Y\in {\cal F}$ the property ${\cal P}$ is satisfied at $x \in Y$ in $X$ if (and only if) the "natural restriction" of ${\cal P}$ on $Y$  is satisfied at $x$. The restriction of ${\cal P}$ on $Y$ is usually done in an evident way. In such a situation we also say that ${\cal P}$ is {\it reducible} by ${\cal F}$. Roughly, separable determinacy reduces the verification of some property in the space $X$ (the latter might not be separable)  to  the verification of the same property on separable subsets of $X$. Examples of such separably reducible properties are some classical properties such as  the reflexivity of Banach spaces, Radon-Nikod\'ym property in Banach spaces, Asplundness of Banach spaces,  Fr\'echet differentiability of continuous functions in Banach spaces (see \cite{Pr} for the latter), Fr\'echet subdifferentiability of functions in Banach spaces (see \cite{FI1}) and many others (see  the book \cite{LPT}).

In this paper we will be interested in separably reducible properties by families which are even richer than merely  covering the underlying spaces. Namely, a family ${\cal R}$ of closed separable subspaces of $X$ is called {\it rich} (see \cite{BoMo}) if the following two properties are true:
\begin{enumerate}
  \item[{\rm (i)}]  for every separable subspace $Y_0$ of $X$ there is an element $Y\in {\cal R}$ such that $Y_0\subset Y$ (cofinality).
  \item[{\rm (ii)}]  if $\{Y_i:i\ge 1\}$ is an increasing sequence of elements of ${\cal R}$  (that is, $Y_i\subset Y_{i+1}$ and $Y_i\in {\cal R}$ for every $i\ge 1$) then $Y=\overline{ \cup_{i\ge 1}Y_i}$ is an element of ${\cal R}$ ($\sigma$-closedness).
\end{enumerate}

The terminology {\it rich} for such families is justified by the following fact:

 \begin{prop}\label{rich_intersection}{\rm  (\cite{BoMo})}
 The intersection of countably many rich families of closed separable subspaces of $X$ is again a rich family of closed separable subspaces of $X$
 \end{prop}

Thus, in particular, if we have two distinct properties (or even  countably many properties!) that are separably reducible by rich families, then the simultaneous satisfaction of all  properties is also separably reducible by the rich family which is the intersection of the corresponding families. Sometimes, when we are in the setting of a Banach space, the (rich) families of closed separable subspaces are constituted by closed separable linear subspaces of $X$.

Results involving separable reducibility by rich families can be found in \cite{BoMo} (for minimality of Clarke  subdifferential and integrability of Lipschitz functions in Banach spaces), \cite{FI2} (for continuity and Fr\'echet subdiffertiability of functions in Banach spaces) in \cite{CF} (for Fr\'echet subdiferrentiability in Asplund spaces) and in \cite{FIR} (for metric regularity) and in many other papers, where the property of separable reducibility is considered. A number  of various examples of separably reducible properties can be found also in the lecture notes \cite{F}.

In this paper we will continue the above lines of research with the investigation of the separable determinacy by rich families of other important properties in metric spaces and Banach spaces. Among them are the local Lipschitz property for functions in such spaces and some functors transforming functions into another ones (like slopes). To do this, in the next section we propose a general scheme of separable reducibility which we apply in the subsequent sections in  concrete situations. The paper is an attempt to propose some unification in these considerations but it is far from being exhaustive concerning the rich  variety of  separably determined properties.

\section{A general aproach for separable determinacy}

Let $(X,d)$ be  a  metric space. As usual, the symbol $B(x,r)$, for $x\in X$ and $r>0$, will stand for the open ball in $(X,d)$ centered at $x$ and with radius $r$. Let $(P,\rho)$ be a given separable metric space and fix a countable dense subset $Q$ of $P$. Put $Z:=X\times P$ and consider on the latter  some usual product metric. In such a setting we will suppose the existence of the following two mappings with the corresponding properties:

\begin{enumerate}

\item[(A1)] There is a set-valued mapping $G:Z\tos X^l$, for some fixed positive integer $l\ge 1$, with the property:

\begin{enumerate}
  \item[(A1.1)] for every $z=(x,p)\in Z$ the set $G(z)$ is nonempty and for any $u\in G(z)$ and any sequence $(x_n)_n\subset X$ which converges to $x$, there is a sequence $(p_n)_n \subset Q$ such that  $p_n\to p$ and $u\in G(z_n)\subset G(z)$ eventually (here $z_n=(x_n,p_n)$ for every $n\ge 1$).
\end{enumerate}
\end{enumerate}

\begin{enumerate}
\item[(A2)] There is a function $\Phi:Z\times X^l\to\Rinf$ such that for every $z\in Z$ we have:
\begin{enumerate}
\item[(A2.1)] for any $u\in G(z)$ if $(z_n)_n\subset Z$ converges to $z$ then $\Phi(z_n,u)\to \Phi(z,u)$ (continuity of $\Phi(\cdot,u)$ at $z$ for every $u\in G(z)$);
\item[(A2.2)] for any nonempty set $A$ which intersects $G(z)$ and for any sequence $(z_n)_n\subset Z$ which converges to $z$ we have
$$\ds
\limsup_{z_n\to z}\,\, \sup\{\Phi(z_n,u): u\in A\cap G(z)\}\le  \sup\{\Phi(z,u): u\in A\cap G(z)\}.
$$
\end{enumerate}
\end{enumerate}

We notice  that if the continuity of $\Phi(\cdot,u)$ at $z\in Z$ from (A2.1) is uniform on $u\in A\cap G(z)$, then the condition (A2.2) is fulfilled.

\medskip

With the above setting we have the following general separable reducibility result:

\begin{theo}\label{general_theorem}
Let $X$, $P$ and  $Q$ be as above. Suppose that (A1) and (A2) are fulfilled with their sub-properties.

Then the family ${\cal R}$ of closed separable subspaces of $X$ such that for every $Y\in {\cal R}$, for every  $x\in Y$ and for every $p\in P$  we have  for $z=(x,p)$ that $Y^l\cap G(z)\neq\emptyset$ and that
$$
\sup\{\Phi(z,u): u\in G(z)\}= \sup\{\Phi(z,u): u\in Y^l\cap  G(z)\}
\leqno{(*)}
$$
is a rich family.

\end{theo}

\proof
We prove first that the family ${\cal R}$ satisfies the condition (i) (cofinality) of the definition of rich family.  To  this end, let us, for any $z=(x,p)\in Z$,  fix a countable set $D(z)\subset  G(z)$ such that

$$
\sup\{\Phi(z,u): u\in G(z)\}=\sup\{\Phi(z,u): u\in D(z)\}.
$$

Let now $Y_0$ be a closed separable subspace of $X$ and $C_0$ be a countable set in $X$ such that $Y_0=\overline{C_0}$. Put
$$
C_1:=C_0\bigcup\Big\{\bigcup_{k=1}^l \pi_k(D(x,p)): x\in C_0, p\in Q\Big\},
$$
where $\pi_k$, $k=1,\ldots,l$, are the canonical projections of $X^l$ onto the $k$-th copy of $X$. Obviously $C_1$ is a countable subset of $Y$  containing $C_0$. Set $Y_1:=\overline{C_1}$. Continuing in this way by induction,  we obtain an increasing sequence of countable sets $C_0\subset C_1\subset\cdots\subset X$ such that for every $i\ge 0$
$$
\ds
C_{i+1}:=C_i\bigcup\Big\{\bigcup_{k=1}^l\pi_k (D(x,p)): x\in C_i, \, \, p \in Q \Big\}.
$$
Let $Y_i:=\overline{C_i}$, $i\ge 0$, and put $Y:=\overline{\cup_{i\ge 0}  Y_i}=\overline{\cup_{i\ge 0}  C_i}$. The latter is a closed separable subspace of $X$ and we will show that $Y\in {\cal R}$. Indeed, let $x\in Y$ and $p\in P$. We have to prove that $Y^l\cap G(z)\neq\emptyset$ and that the left hand side of (*) is less or equal to the right hand one.

To prove this fix an element $u\in G(z)$ where $z=(x,p)$.  Let $(x_n)_n\subset \cup_{i\ge 0} C_i$ be a sequence which converges to $x$. By condition (A1.1) there is a sequence $(p_n)_n \subset Q$ which converges to $p$ and we have that for the sequence $z_n=(x_n,p_n)$, $n\ge 1$,
$$
u\in G(z_n)\subset G(z) \mbox{ eventually. }
\leqno{(2.1)}
$$
On the other hand, let $i_n$ be such that $x_n\in C_{i_n}$, $n\ge 0$. Then, by the definition of $C_{i_n+1}$, we have that  $\cup_{k=1}^l\pi_k (D(x_n,p_n))$ belongs to $C_{i_n+1}$ and thus, $D(x_n,p_n)$ is contained in $Y^l$ for each $n\ge 1$. By (2.1), this, in particular, implies  that $Y^l$ intersects $G(z)$. We proved that $Y^l\cap G(z)$ is nonempty.

Further, suppose first  that $\Phi(z,u)$ is finite. Fix any $\eps>0$. Let $n$ be so large such that (2.1) is fulfilled and also,  using  (A2.1) and (A2.2), that for any such $n$ we have in addition,
$$
\Phi(z,u)\le\Phi(z_n,u)+\eps
\leqno{(2.2)}
$$
 and
$$
\sup\{\Phi(z_n,v): v\in Y^l\cap G(z)\}\le \sup\{\Phi(z,v): v\in Y^l\cap G(z)\}+\eps.
\leqno{(2.3)}
$$

Then, using the last two inequalities,  the fact that $D(x_n,p_n)\subset Y^l$ for every $n$, the definition of $D(x_n,p_n)$ and that by (2.1) for $n$ large enough we have $D(x_n,p_n)\subset  G(z)$, we obtain that the following chain of inequalities holds true for large $n$:

$$
\begin{array}{ll}
\vspace{6pt}
\Phi(z,u)&\le \Phi(z_n,u)+\eps\le \sup\{\Phi(z_n,v): v\in G(z_n)\}+\eps\\

\vspace*{6pt}
 &=\sup\{\Phi(z_n,v): v\in D(z_n)\}+\eps \le \sup \{\Phi(z_n,v): v\in Y^l\cap G(z)\} +\eps \\
 \vspace{6pt}
& \le \sup \{\Phi(z,v): v\in Y^l\cap G(z)\} +2\eps.
\end{array}
$$
Since $\eps>0$ was arbitrary this shows that $\Phi(z,u)$ is less or equal to the right hand side of (*). In the case when $\Phi(z,u)=+\infty$ the reasoning is similar: For any $M>0$ we have that $\Phi(z_n,u)>M$ eventually and then we use the same argument (and chain of inequalities) as above to show that $\sup \{\Phi(z,v): v\in Y^l\cap G(z)\}>M$, which, since $M>0$ was arbitrary,  implies the equality in (*).

\smallskip

To prove the second property (ii) ($\sigma$-closedness) from the definition of rich family suppose that $Y_1\subset Y_2\subset\cdots\subset X$ is an increasing sequence  of elements from ${\cal R}$. Put $Y:=\overline{\cup_{i\ge1} Y_i}$, which obviously is a closed separable subset of $X$  and let $x\in Y$ and $p\in P$. Let $u\in G(z)$ (with $z=(x,p)$) be arbitrary.  Let $(x_n)_n\subset \cup_{i\ge 0} Y_i$ be a sequence which converges to $x$. Condition (A1.1) gives the existence of a sequence $(p_n)_n \subset Q$ which converges to $p$ and we have that the sequence $z_n=(x_n,p_n)$, $n\ge 1$, satisfies (2.1). Suppose that $x_n\in Y_{i_n}$ for some $i_n$, $n\ge 1$. Since $Y_{i_n}^l\cap G(z_n)\neq\emptyset$ and $G(z_n)\subset G(z)$ eventually, we obtain that $Y^l\cap G(z)\neq\emptyset$.

Again, suppose first that $\Phi(z,u)$ is finite and fix any  $\eps>0$. Let $n$ be so large that, in addition to (2.1) we have also (2.2) and (2.3) fulfilled.

Then,  having in mind that $Y_i\in {\cal R}$ for any $i\ge 1$, we obtain the following chain of inequalities for large $n$:
$$
\begin{array}{ll}
\vspace{6pt}
\Phi(z,u)&\le \Phi(z_n,u)+\eps\le \sup\{\Phi(z_n,v): v\in G(z_n)\}+\eps\\
\vspace{6pt}
 &=\sup\{\Phi(z_n,v): v\in Y_{i_n}^l \cap G(z_n)\}+\eps \le \sup \{\Phi(z_n,v): v\in Y^l\cap G(z)\} +\eps \\
 \vspace{6pt}
& \le \sup \{\Phi(z,v): v\in Y^l\cap G(z)\} +2\eps
\end{array}
$$
obtaining the desired inequality in (*) since $\eps>0$ was arbitrary. The case when $\Phi(z,u)$ is infinite is treated analogously. The proof is completed.

\eproof

\begin{rem}{\rm
Observe that, if we are in the setting of a Banach space $X$ and we want to work only with closed separable linear subspaces of $X$, then in the above constructions we simply have to consider instead of $C_i$ the sets $C_i'$ which consist of all rational linear combinations of elements from $C_i$.
}
\end{rem}

We  have a "dual version" of the above result replacing the supremum with infimum provided we suppose instead of (A2.2) the following property:

\begin{enumerate}
\item[(A2.2$'$)] for any nonempty set $A$ which intersects $G(z)$ and for any sequence $(z_n)_n\subset Z$ which converges to $z$ we have
$$\ds
\liminf_{z_n\to z}\,\, \inf\{\Phi(z_n,u): u\in A\cap G(z)\}\ge  \inf\{\Phi(z,u): u\in A\cap G(z)\}.
$$
\end{enumerate}

Then, with  similar arguments we can prove the following  general separable reducibility result:

\begin{theo} \label{general_theorem_bis}
Let $X$, $P$ and  $Q$ be as above. Suppose that (A1) and (A2) are fulfilled with their sub-properties replacing (A2.2) with (A2.2$'$) and supposing that $\Phi$ is with values in $\R\cup\{-\infty\}$.  Then the family ${\cal R}$ of closed separable subspaces of $X$ such that for every $Y\in {\cal R}$, for every  $x\in Y$ and for every $p\in P$ we have for $z=(x,p)$ that $Y^l\cap  G(z)\neq\emptyset$ and that
$$
\inf\{\Phi(z,u): u\in G(z)\}= \inf\{\Phi(z,u): u\in Y^l\cap  G(z)\}
\leqno{(**)}
$$
is a rich family
\end{theo}

At the end of this section we will investigate separable reducibility in the above general form when the space $X$ is of  the form of Cartesian product  $X=X_1\times X_2$, where $(X_1,d_1)$ and $(X_2,d_2)$ are two metric spaces and on $X$ we consider some usual product metric generated by $d_1$ and $d_2$. In particular, we will see that in such a case  we can restrict ourselves to take  rich families  consisting only of elements of the form $Y_1\times Y_2$, where $Y_i\in {\cal S}(X_i)$, $i=1,2$. Namely, we have the following result:

\begin{theo}\label{general_theorem_product}
Let $P$ and  $Q$ be as above and let $X=X_1\times X_2$, where $(X_1,d_1)$ and $(X_2,d_2)$ are metric spaces. Suppose that $G:X_1\times P\tos X_1^l$ is a mapping satisfying (A1) with respect to  $Z:=X_1\times P$. Let  $\Phi: Z\times X_1^l\times X_2\to\Rinf$ be a function such that for every $y\in X_2$ the function $\Phi(\cdot, \cdot, y)$ satisfies (A2.1) and (A2.2) with respect $Z$, and in addition for every $z\in Z$ the function $\Phi(z,u,\cdot)$ is continuous at any $y\in X_2$ uniformly on $u\in G(z)$.   Then, the family ${\cal R}$ of closed separable subspaces of $X$ such that every $Y\in {\cal R}$ is of the form $Y=Y_1\times Y_2$, where $Y_i\in {\cal S}(X_i)$, $i=1,2$, and, in addition,  for every  $(x,y)\in Y=Y_1\times Y_2$  and for every $p\in P$  we have  for $z=(x,p)$ that $Y_1^l\cap G(z)\neq\emptyset$ and that
$$
\sup\{\Phi(z,u,y): u\in G(z)\}= \sup\{\Phi(z,u,y): u\in Y_1^l\cap  G(z)\}
\leqno{(***)}
$$
is a rich family.
\end{theo}

\proof

To prove the cofinality of ${\cal R}$,  let $Y_0$ be a (closed) separable subspace of $X=X_1\times X_2$. Then, we can find closed separable subspaces $Y_1^0\in {\cal S}(X_1)$ and $Y_2\in {\cal S}(X_2)$ such that $Y_0\subset Y_1^0\times Y_2$. Take a countable set $C\subset Y_2$ which is dense in $Y_2$. For every fixed $y\in C$, let ${\cal R}_y$ be the rich family of closed separable subsets of $X_1$ which is provided by Theorem \ref{general_theorem} for the function $\Phi(\cdot,\cdot,y)$ and which for every $Y_1\in {\cal R}_y$, every  $x\in Y_1$ and $p\in P$  will satisfy $Y_1^l\cap G(z)\neq\emptyset$ and  (***) for the fixed $y$ and $z=(x,p)$. Set ${\cal R}':=\cap\{{\cal R}_y: y\in C\}$. This is a rich family of closed separable spaces in $X_1$ according to Proposition \ref{rich_intersection}. Find an element $Y_1\in {\cal R}'$ such that $Y_1^0\subset Y_1$. It can be seen, using the fact that for every $z=(x,p)\in Z$ the function $\Phi(z,u,\cdot)$ is continuous at any $y\in X_2$ uniformly on $u\in G(z)$ that the closed separable subspace $Y_1\times Y_2$ belongs to ${\cal R}$.

Let us further prove the $\sigma$-completeness of ${\cal R}$. To this end, let $L_1\subset \cdots L_n\subset \cdots$ and $M_1\subset\cdots \subset M_n\subset \cdots$ be increasing sequences of closed separable subsets of $X_1$ and $X_2$ correspondingly, such that $L_n\times M_n\in {\cal R}$ for every $n\ge 1$. Set $Y_1:=\overline{\cup_{n\ge 1} L_n}$ and $Y_2:=\overline{\cup_{n\ge 1} M_n}$. We will prove that $Y_1\times Y_2\in {\cal R}$. Indeed, let $C$ be a dense countable set in $Y_2$. Fix $y\in C$ and let again, ${\cal R}_y$ be the rich family provided by Theorem \ref{general_theorem} for the function $\Phi(\cdot,\cdot,y)$ and which satisfies that for every $Y_1\in {\cal R}_y$ and for every  $z=(x,p)\in Y_1\times P$ the set $Y_1$ intersects $G(z)$ and (***) is true. Then, obviously $L_n\in {\cal R}_y$ for every $n\ge 1$ and, therefore, $Y_1\in {\cal R}_y$ as well, because ${\cal R}_y$ is rich. Now, the property that for every $z=(x,p)\in Z$ the function $\Phi(z,u,\cdot)$ is continuous at each $y\in X_2$ uniformly on $u\in G(z)$ implies that we have (***) for each $y\in Y_2$ and therefore,  $Y_1\times Y_2$ belongs to ${\cal R}$. The proof is completed.

\eproof

\section{Separable reducibility of the continuity and the Lipschitz property}

In this section we will investigate the separable determiness of the continuity and the Lipschitz property of  functions  in metric spaces applying our general scheme. As we mentioned, the separable reducibility of the continuity property can be found, say in \cite{FI2}. Here is how it can be derived from the above general theorems.

\begin{theo}
Let $f:(X,d)\to \R$ be a function defined in the metric space $X$. Then there is a rich family ${\cal R}$ of closed separable subsets of $X$ such that for every $Y\in {\cal R}$ and  every $x\in Y$ we have
\begin{enumerate}
\item[{\rm (a)}] $\liminf_{y\to x}f(y)=\liminf_{Y\ni y\to x}f(y)$;
\item[{\rm (b)}] $\limsup_{y\to x}f(y)=\limsup_{Y\ni y\to x}f(y)$;
\item[{\rm (c)}] $\lim_{y\to x}f(y)=f(x)$ if, and only if, $\lim_{Y\ni y\to x}f(y)=f(x)$.
\end{enumerate}
\end{theo}

\proof
To prove (a) we apply Theorem \ref{general_theorem_bis}. Let  $P:=\R_{++}:=\{r\in \R: r>0\}$ and  $Q$ be the set of all strictly positive rational numbers. For $Z:=X\times P$ and $z=(x,r)\in Z$, let $G(z):=B(x,r)$. The mapping $G$ is easily seen to satisfy (A1). Let  the function $\Phi:Z\times X\to \R$ be defined as $\Phi(z,u):=f(u)$, $z=(x,r)\in Z$, $u\in X$. The verification of (A2.1) and A(2.2$'$) are immediate. Therefore, there is a rich family ${\cal R}_1$ of closed separable subsets of $X$ such that for any $Y\in {\cal R}_1$, for any $x\in Y$ and any $r>0$ we have $\inf\{f(u): u\in B(x,r)\}= \inf\{f(u): u\in Y\cap B(x,r)\}$. By the definition of $\liminf$ it follows that for every  $Y\in {\cal R}_1$ and  for any $x\in Y$ we have $\liminf_{y\to x}f(y)=\liminf_{Y\ni y\to x}f(y)$.
In a similar way, using Theorem \ref{general_theorem} we obtain a rich family ${\cal R}_2$ of closed separable subsets of $X$ such that for every  $Y\in {\cal R}_2$ and  for any $x\in Y$ we have $\limsup_{y\to x}f(y)=\limsup_{Y\ni y\to x}f(y)$. To complete the proof, it remains to take the intersection of the two rich families ${\cal R}:= {\cal R}_1\cap {\cal R}_2$ which, according to Proposition \ref{rich_intersection}, is again a rich family.
\eproof

We continue now with the investigation of the separable determinacy of the Lipschitz property.  We start with the following proposition:

\begin{prop}
Let $f:(X,d)\to \R$ be a function defined in the metric space $X$. Then there is a rich family ${\cal R}$ of closed separable subsets of $X$ such that for every $Y\in {\cal R}$, every $x\in Y$ and any $r>0$ we have
$$\ds
\begin{array}{ll}\ds
\sup\Bigg\{\frac{|f(u_1)-f(u_2)|}{d(u_1,u_2)}: u_1\neq u_2, u_1,u_2\in B(x,r)\Bigg\}=&\hspace*{4cm}\\
\ds
\qquad \qquad \qquad \qquad=\sup\Bigg\{\frac{|f(u_1)-f(u_2)|}{d(u_1,u_2)}: u_1\neq u_2, u_1,u_2\in Y\cap B(x,r)\Bigg\}.&\\
\end{array}
$$

\end{prop}

\proof

The proof is a straightforward application of our Theorem 2.1. Put $P:=\R_{++}=\{r\in \R: r>0\}$ with the usual metric in $\R$ and let $Q$ be the set of all strictly positive rational numbers.  Let us for $z=(x,r)\in Z:=X\times P$ define $G(z):=\{(u_1,u_2)\in X^2: u_1\neq u_2, u_1,u_2\in B(x,r)\}$. It is easy to see that $G$ satisfies (A1.1). Let $\Phi:Z\times X^2\to \R$ be defined by $\Phi((x,r),u):=|f(u_1)-f(u_2)|/d(u_1,u_2)$, where $z=(x,r)\in Z$ and $u=(u_1,u_2)\in X^2$. Since the values of $\Phi$ do not depend on $z$, the conditions (A2.1) and (A2.2) are trivially fulfilled. Thus, it remains to apply Theorem 2.1
\eproof

With the above proposition in hand, we can prove that the (local) Lipschitzian property of a function in $X$ is separably reducible. Using (with  simplified notation) the above used function $\Phi(u_1,u_2):=|f(u_1)-f(u_2)|/d(u_1,u_2)$, $u_1,u_2\in X$, recall that for a function $f:(X,d)\to\R$ the {\it Lipschitz modulus} $\Lip f(x)$ of $f$ at $x\in X$ is defined as
$$
\Lip f(x):=\inf\Bigg\{c>0: \exists r>0:\,  \sup\Big\{\Phi(u_1,u_2): u_1\neq u_2\in B(x,r)\Big\}\le c \Bigg\}.
$$

The following theorem shows that the local Lipschitz property of a function in metric spaces is separably reducible. As usual, for a function $f:X\to\R$ and a set $A\subset X$, the symbol $f|_A$  means the restriction of $f$ on $A$.

\begin{theo}
Let $f:(X,d)\to \R$ be a function defined in $X$. Then there is a rich family ${\cal R}$ of closed separable subsets of $X$ such that for every $Y\in {\cal R}$ and  every $x\in Y$ we have
$$
\Lip f(x)=\Lip (f|_Y)(x).
$$
\end{theo}

\proof
Given $f$ let ${\cal R}$ be the rich family provided by Proposition 3.2. Let $Y\in {\cal R}$ and $x\in Y$ be fixed.  Let $\Phi(\cdot,\cdot)$ be as above and consider the two sets of possible Lipschitz constants on balls (centered at $x$)  in $X$ and in  $Y$ respectively:
$$
A:=\Big\{ c>0: \exists r>0:  \sup\{\Phi(u_1,u_2): u_1\neq u_2\in B(x,r)\}\le c\Big\}
$$
and
$$
A_Y:=\Big\{ c>0: \exists r>0:  \sup\{\Phi(u_1,u_2): u_1\neq u_2\in Y\cap B(x,r)\}\le c\Big\}.
$$

Obviously, $A\subset A_Y$ and therefore, $\inf A\ge \inf A_Y$ (with the convention $\inf \emptyset = +\infty$). To prove  the inverse inequality, take and fix some $\eps>0$ and suppose $A_Y\neq\emptyset$ (otherwise we are done). Let $\bar c \in A_Y$ be such that $\bar c<\inf A_Y +\eps$. Therefore, for some $r>0$ we have
$$
\sup\{\Phi(u_1,u_2): u_1\neq u_2\in Y\cap B(x,r)\}\le \bar c.
$$
By Proposition 3.2 we have
$$
\sup\{\Phi(u_1,u_2): u_1\neq u_2\in B(x,r)\}=\sup\{\Phi(u_1,u_2): u_1\neq u_2\in Y\cap B(x,r)\}
$$
and therefore, $\bar c\in A$. Consequently, $\inf A\le \bar c<\inf A_Y+\eps$ and since $\eps>0$ was arbitrary we obtain that $\inf A\le \inf A_Y$, which completes the proof of the theorem.

\eproof

\begin{rem} {\rm
If $X$ is a Banach space, then, having in mind Remark 2.2, the seprable determinacy of the Lipschitz property can be obtained only by closed separable linear subspaces of $X$.

}
\end{rem}

\section{Separable reducibility of slopes}

Let $(X,d)$ be a metric space and $f:X\to\Rinf$ be a proper extended real-valued function in it. Recall that the {\it slope} of $f$ at $x\in X$ is given by the expression:
$$
\ds
|\nabla f|(x):=\limsup_{x\neq u\to x}\frac{(f(x)-f(u))^+}{d(x,u)}, \qquad x\in X,
$$
where the symbol $s^+$ for a real number $s$ has the usual meaning (0 of $s\le 0$ and $s$ otherwise) and we adopt the convention $+\infty-\infty=0$. We will show that the slope of a function is a functor which is separably reducible by rich families of closed separable sets in $X$. By this we simply mean that if $Y$ is a member of such a family and $x\in Y$ then $|\nabla f|(x)=|\nabla (f|_Y)|(x)$. To prove it, we need the following proposition which is a consequence of our main theorem. Below, given  a point $x\in X$ and real positive numbers $r,s$ such that $0<r<s$ we denote by $T(x,r,s):=\{u\in X: r<d(x,u)<s\}$ the {\it torus} in $X$ generated by these parameters.

\begin{prop}
Let $f:X\to\Rinf$ be  a proper extended-real valued function in $X$. Then the family ${\cal R}$ of closed separable sets in $X$ such that for every $Y\in {\cal R}$, for every $x\in Y$, for every $t\in \R$ and every $r,s$ with $0<r<s$ we have $Y\cap T(x,r,s)\neq\emptyset$ and
$$
\sup\Bigg\{\frac{(t-f(u))^+}{d(x,u)}: u\in T(x,r,s)\Bigg\}=\sup\Bigg\{\frac{(t-f(u))^+}{d(x,u)}: u\in Y\cap T(x,r,s)\Bigg\}
$$
is a rich family.
\end{prop}

\proof
To apply our general  theorem, let $P:=\{(t,r,s)\in \R^3: 0<r<s\}$ with the inherited standard metric from $\R^3$. Consider as $Q$ the triples $(t,r,s)$  belonging to $P$ whose coordinates are rational numbers.  Let $Z:=X\times P$ and for $Z\ni z=(x,p)\in X\times P$, $p=(t,r,s)$, define $G:Z\tos X$ by $G(z):=T(x,r,s)$. It is easily seen that the so defined $G$ satisfies (A1.1). Define now the function $\Phi:Z\times X\to \Rinf$ by
$$
\Phi (z,u):= \frac{(t-f(u))^+}{d(x,u)}, \qquad z=(x,(t,r,s))\in Z, u\in X.
$$
Observe that for $z=(x,(t,r,s))\in Z$ and $u\in G(z)$ the distance $d(x,u)$ is always greater than $r>0$. Therefore, for any $u\in G(z)$ the function $\Phi(\cdot,u)$ is continuous at $z\in Z$, thus (A2.1) is fulfilled. To check (A2.2), take $z\in Z$ and let $A$ be a subset of $X$ intersecting $G(z)$. If $\sup\{\Phi(z,u):u\in A\cap G(z)\}=+\infty$ then (A2.2) is obviously satisfied. If the latter supremum is finite then this implies, in particular, that $f$ is bounded below on $A\cap G(z)$ (remember that $0<r<d(x,u)<s$ for any $u\in G(z)$). In the latter case it  can be seen that $\Phi(\cdot,u)$ is continuous in $z$ uniformly on $u\in A\cap G(z)$. But in this situation, as it was mentioned  after the formulation of (A2.2), the latter property is satisfied. Therefore, we can apply our Theorem 2.1 to get the conclusion of the proposition, which completes the proof.

\eproof

Having in mind the above, we can prove the separable determiness of slopes.

\begin{theo}
Let $(X,d)$ be a metric space and $f:X\to\Rinf$ be  a proper extended real-valued function in $X$. Then there is a rich  family ${\cal R}$ of closed separable sets in $X$ such that for every $Y\in {\cal R}$ and  for every $x\in Y$ we have
$$
|\nabla f|(x)=|\nabla (f|_Y)|(x).
$$
\end{theo}

\proof
Let ${\cal R}$ be the rich family provided by the previous proposition. Then for every $Y\in {\cal R}$ and $x\in Y$ we have
$$
\begin{array}{ll}
|\nabla f|(x)& \ds= \inf_{s>0} \sup_{r\in (0,s)} \sup \Bigg\{\frac{(f(x)-f(u))^+}{d(x,u)}: u\in T(x,r,s)\Bigg\}\\
 &\ds =\inf_{s>0} \sup_{r\in (0,s)} \sup \Bigg\{\frac{(f(x)-f(u))^+}{d(x,u)}: u\in Y \cap T(x,r,s)\Bigg\}=|\nabla (f|_Y)|(x),\\
\end{array}
$$
Indeed, when $x\in \dom f$ this is immediate from the previous proposition. To see that it is true also for  $x\notin \dom f$, let us just observe that the previous proposition shows also  that the case $\dom f\cap T(x,r,s)\neq \emptyset$ and $\dom f\cap Y\cap  T(x,r,s)=\emptyset$ is not possible and thus, we have that both  inner supremums above are either 0 or $+\infty$.

\eproof

We continue with the study of the separable reducibility of slopes in the following setting. Let $(X_1,d_1)$ and $(X_2,d_2)$ be two metric spaces and  $f:X_1\times X_2\to \R$ be an extended real-valued function on the Cartesian product of these two spaces, which is Lipschitz in the second variable, in a sense that there is $k>0$ such that
$$
|f(x,y')-f(x,y'')|\le k \, \, d_2(y',y''), \qquad \mbox{for every } x\in X_1.
$$
On the Cartesian product $X_1\times X_2$ we consider some usual product metric generated by $d_1$ and $d_2$. Then, we have the following theorem:

\begin{theo}
Let $(X_1,d_1)$, $(X_2,d_2)$ and $f:X_1\times X_2\to \R$ be as above. Then there is a rich family ${\cal R}$ of  separable subspaces of $X_1\times X_2$ of the form $Y_1\times Y_2$, $Y_i\in {\cal S}(X_i)$, $i=1,2$, such that for every $Y_1\times Y_2\in {\cal R}$, and for every $(x,y)\in Y_1\times Y_2$ we have $|\nabla f(\cdot,y)|(x)=|\nabla (f|_{Y_1\times Y_2})(\cdot,y)|(x)$.
\end{theo}

\proof
Let $d$ be some product metric in $X:=X_1\times X_2$. We proceed first as in the proof of Proposition 4.1:  let $P:=\{(t,r,s)\in \R^3: 0<r<s\}$ and $Q$ be the set of  the rational triples $(t,r,s)$  belonging to $P$.  Let $Z:=X_1\times P$ and for $Z\ni z=(x,p)\in X_1\times P$, $p=(t,r,s)$, define $G:Z\tos X_1$ by $G(z):=T(x,r,s)$, the latter being the corresponding torus in $X_1$. The mapping  $G$ satisfies (A1.1). Define the function $\Phi:Z\times X_1\times X_2\to \Rinf$ by
$$
\Phi (z,u,y):= \frac{(t-f(u,y))^+}{d_1(x,u)}, \qquad z=(x,(t,r,s))\in Z, \, u\in X_1, \,  y\in X_2.
$$
It can be easily seen, using the Lipschitz property of $f$ wrt to the second variable,  that for every $z\in Z$ the function $\Phi(z,u,\cdot)$ is continuous at each $y\in X_2$ uniformly on $u\in G(z)$. As in the proof of Proposition 4.1, for every fixed $y\in Y_2$, the function $\Phi(\cdot,\cdot,y)$ satisfies (A2.1) and (A2.2). Then, by Theorem \ref{general_theorem_product}, there is a rich family ${\cal R}$ of separable subspaces of $X$ of the form $Y_1\times Y_2$, $Y_i\in{\cal S}(X_i)$, $i=1,2$,  such that for every $Y_1\times Y_2\in{\cal R}$, for any $(x,y)\in Y_1\times Y_2$ and for any $(t,r,s)\in P$ we have
$$
\sup\Bigg\{\frac{(t-f(u,y)^+}{d_1(x,u)}: u\in T(x,r,s)\Bigg\}=\sup\Bigg\{\frac{(t-f(u,y))^+}{d_1(x,u)}: u\in Y_1\cap T(x,r,s)\Bigg\}.
$$

Further, the proof follows the steps from the proof of the previous theorem.
\eproof

\bigskip

\noindent
M.Fabian\\
Institute of Mathematics, Czech Academy of Sciences, \v Zitn\'a 25, 115 67 Prague, Czech Republic\\
e-mail: fabian@math.cas.cz

\bigskip

\noindent
A. Ioffe\\
Department of Mathematics, Technion-Israel Institute of Technology, Haifa 32000, Israel\\
e-mail: alexander.ioffe38@gmail.com

\bigskip

\noindent
J.P. Revalski\\
Institute of Mathematics and Informatics, Bulgarian Academty of Sciences, Acad. G. Bonchev str., block 8, 1113 Sofia, Bulgaria\\
\noindent
and\\
\noindent
Bulgarian Academy of Sciences, 1, 15th November str., 1040 Sofia, Bulgaria\\
e-mail: revalski@math.bas.bg

\end{document}